\documentclass[reqno]{amsart}
\newif\ifpdf
\ifx \pdfoutput \undefined
  \usepackage{graphicx}
  \pdffalse
\else
  \usepackage[pdftex]{graphicx}
  \DeclareGraphicsExtensions{.pdf,.png,.jpg,.jpeg}
\fi
\usepackage{geompsfi}
\usepackage{amsmath,amsfonts,amssymb}
\usepackage[english]{babel}
\usepackage{natbib}
\let\cite\citep
\bibpunct();a{},

\def\Omega{V}           
\def\dOmega{S}          
\def\avj{\alpha_{v_j}}      
\def\kk{\kappa_k}
\def\kp{\kappa_p}
\def\ak{\alpha_k}
\def\ap{\alpha_p}
\def\afi{\alpha_{f_i}}
\def\ad{\alpha_D^{}}
\def\adi{\alpha_{D_i}^{}}
\def\at{\alpha_t}
\def\al{\alpha_L}
\def\r{v}               

\def\X{Y}
\def\XP{{Y\!P}}
\def\Tend{\tau_{\mathrm{end}}}
\def\tend{t_{\mathrm{end}}}
\long\def\drop#1{}
\newcounter{remark}[section]
\def\theremark{\thesection.\arabic{remark}}
    {\par\medbreak\refstepcounter{remark}\textsc{Remark \theremark~}}%
    {\par\medskip}%
\def\pref#1{\ref{#1}}
\def\R{{\mathbb R}}
\def\mod#1{\left|#1\right|}

\def\ld{,\ldots,}
\def\equil#1{%
  \quad
  \mathop{\rightleftarrows}\limits^{{#1}}%
  \quad
}
\begin{document}
\title[Control of Heterogeneous Networks]{Control of Spatially  
Heterogeneous and Time-Varying 
Cellular Reaction Networks: A~New Summation Law}
\author{Mark A. Peletier, Hans V. Westerhoff, and Boris N. Kholodenko}
\thanks{\noindent Centrum voor Wiskunde en Informatica, Amsterdam, The Netherlands; \texttt{peletier@cwi.nl}\\
Vrije Universiteit Amsterdam, The Netherlands; \texttt{hw@bio.vu.nl}\\
Thomas Jefferson University, Philadelphia, USA; \texttt{Boris.Kholodenko@mail.tju.ed}}

\begin{abstract}

A hallmark of a plethora of intracellular signaling pathways is the
spatial separation of activation and deactivation processes that
potentially results in precipitous gradients of activated
proteins. The classical Metabolic Control Analysis (MCA), which
quantifies the influence of an individual process on a system variable
as the control coefficient, cannot be applied to spatially separated
protein networks. The present paper unravels the principles that
govern the control over the fluxes and intermediate concentrations in
spatially heterogeneous reaction networks.  Our main results are two
types of the control summation theorems. The first type is a
non-trivial generalization of the classical theorems to systems with
spatially and temporally varying concentrations. In this
generalization, the process of diffusion, which enters as the result
of spatial concentration gradients, plays a role similar to other
processes such as chemical reactions and membrane transport. The
second summation theorem is completely novel. It states that the
control by the membrane transport, the diffusion control coefficient
multiplied by two, and a newly introduced control coefficient
associated with changes in the spatial size of a system (e.g., cell),
all add up to one and zero for the control over flux and
concentration. Using a simple
example of a kinase/phosphatase system in a spherical cell, 
we speculate that unless active
mechanisms of intracellular transport are involved, the threshold cell
size is limited by the diffusion control, when it is beginning to
exceed the spatial control coefficient significantly.

\drop{
The original focus of Metabolic Control Analysis was on time-independent,
spatially homogeneous metabolic systems. After extensions to 
time-varying systems~\cite{AcerenzaSauroKacser89} and into
the frequency domain~\cite{ReijengaWesterhoffKholodenkoSnoep02}
we here extend the classical MCA results to spatially
heterogeneous systems. 

Our main results are two types of summation theorems. The first type
is a non-trivial but unsurprising 
generalization of the classical flux/concentration summation theorem
to systems with spatially (and temporally) varying concentrations.
In this theorem the process of diffusion, which enters as the result of spatial
concentration gradients, plays a role similar to other
processes such as chemical reactions and membrane transport.

The second type of summation theorem is completely novel, and introduces
a new control coefficient $C_L$ associated with changes in the size of
the system:
\[
2C_{\textrm{diffusion}} + C_{\textrm{membrane transport}}
  + C_L = a,
\]
where $a=1$ for flux control and $a=0$ for concentration control. 
In contrast with the classical result reactive processes do not enter
in this theorem,
and the diffusion control coefficient is multiplied by a factor $2$.

In this paper we derive these two types of summation theorem,
comment on the interpretation of the theorems and the novel
control coefficient $C_L$, and examine in detail a simple example of
a kinase/phosphatase system. 
}
\end{abstract}
\maketitle
\section{Introduction}

Extracellular information received by plasma membrane receptors is
processed, encoded, and transferred to the nucleus through activation
and spatial relocation of multiple signalling components. Receptor
activation triggers signalling responses associated with the
mobilization of a plethora of adapter and target proteins to the 
plasma membrane~\cite{HaughLauffenburger97,KholodenkoHoekWesterhoffTCB}. 
Proteins activated e.g.\ by phosphorylation at the cell
surface travel to stimulate critical regulatory targets at various
cellular sites including the nucleus. During and after the travel,
these proteins are inactivated, e.g.\ through dephosphorylation. 
For instance, a protein
phosphorylated by a plasma-membrane associated protein kinase 
can be dephosphorylated by
a phosphatase in the cytosol or nucleus.  The
transport between cellular locations, where the activation and
inactivation occur, is passive (thermal diffusion), driven by spatial
gradients in the concentrations. 
In some cases, not to be discussed here, motor proteins or 
endocytosis
may be involved~\cite{Kholodenko02}. Previous work estimated that these
spatial gradients may be large,  and that therefore diffusion may
contribute to the control of signal 
transduction~\cite{BrownKholodenko99,KholodenkoBrownHoek00}.

For spatially homogeneous reaction networks, the control over
fluxes and intermediate concentrations has been studied
both experimentally and theoretically (reviewed in~\cite{Fell97}). The
control is generally quantified as the extent to which any type of molecular
process influences a system variable, such as the flux or
concentration. The control coefficients are defined as
the ratios of the fractional changes in the system variable to that
of the biochemical activity which caused the
system change. Mathematically, the stability of a system steady state
is required and the changes are considered as infinitesimally small,
i.e., the coefficients are expressed as the log to log derivatives~\cite{Fell97}.
Important principles underlying the control of
biochemical reaction networks in well-stirred reactors have been
worked out recently. For instance, it can be shown that a large
increase in the activity of a single enzyme does not result in a
substantial increase in the flux for almost any metabolic network;
marked flux increase can be achieved by a concerted modulation of
several pathway reactions~\cite{KacserAcerenza93}. 
In metabolic control analysis (MCA),
this result is related to the so-called summation theorem, which
states that the sum of the enzyme control coefficients adds up to
$1$~\cite{KacserBurns73,HeinrichRapoport73,WesterhoffVanDam87}. 
There is no rate-limiting enzyme, the control is shared between all
network processes.

The present paper unravels the principles that govern the control
pattern in spatially heterogeneous cellular networks. We demonstrate
new properties of control that result from the spatial
aspect of diffusion and separation of signalling reactions. The
control summation theorems, relevant for these networks, are formulated
and proven.

\section{Methods}

\subsection{Description of the system}

We consider living metabolic systems consisting of $n$ internal 
chemical species $Y_1\ld Y_n$,
with concentrations $c_1\ld c_n$. These species interact
via $m$  reactions which can be represented as 
\begin{equation}
\label{reaction:general}
\sum_{i=1}^n N_{ij}^+ Y_i \equil{} \sum_{i=1}^n N_{ij}^- Y_i,
  \qquad j = 1\ld m
\end{equation}
where $N_{ij}^{\pm}$ are the dimensionless forward and reverse stoichiometric
coefficients for the consumption of species $Y_i$ in the forward
and reverse reactions $j$. We write $N_{ij} = N_{ij}^- - N_{ij}^+$
for the net stoichiometric component. Reaction $j$ has rate
$\avj\r_j(c_1\ld c_n)$, so that in the absence of diffusion
the time evolution of the concentrations
$c_i$ should be given by
\[
\frac{\partial}{\partial t} c_i = \sum_{j=1}^m N_{ij} \avj\r_j, \qquad
   i = 1\ld n.
\]
$c_i$ has the dimension of concentration [M], or mole number per unit
volume. 
The rate of the reaction~\pref{reaction:general}, i.e.\ $\r_j$, has the
dimension of concentration per unit time [M/s].  It is a function of
the kinetic parameters of the enzyme catalyzing that reaction (or if
non-enzymatic of the chemical reaction itself), 
e.g.\ $K_M$,
the concentration of the enzyme, and the concentrations of substrates,
products and allosteric modifiers of the reaction. 
The dimensionless parameter $\avj$ is the tool we shall use
to modulate the activity of reaction $\r_j$.

In this paper we will be interested in cases in which concentrations
can \emph{not} be assumed homogeneous in space inside the systems. 
For thermodynamic reasons, living systems need to be open for at least
some chemical compounds. These substances are transported across
the system membrane, usually by transport proteins, i.e.\ effectively 
there is a sink or source for some of the species
on the system boundary, e.g.\ the (plasma) membrane. 
The spatial separation between the sink or source and
the chemistry that takes place in the bulk phase
then creates spatial variations in concentration.

When the system is not
homogeneous and not at steady state, the concentrations need to be
specified as functions of time and space (i.e.\ $c_i(x,t)$), and
boundary fluxes are required.
The principal difficulty
accompanying the spatial inhomogeneity
of the concentrations of any component $Y_i$, which
makes the classical MCA inapplicable, is that the number of `explicit'
variables associated with $Y_i$ becomes infinitely large.
The reactions are
assumed to take place in an enclosed three-dimensional space  that is
represented by a bounded set of $x_1,x_2,x_3$-coordinates, i.e.\ the volume
$\Omega\subset \R^3$. For instance, $\Omega$ can
correspond to the cytoplasm or the mitochondrial matrix. Inside $\Omega$
the species react with each other, according to Eq.~\pref{reaction:general};
in addition they undergo diffusive transport. The concentrations 
$c_i$ of the species are functions of time ($0 < t < \tend$) and space
($x\in\Omega$), and their evolution is governed by the balance
(reaction-diffusion)
equations~\cite{KatchalskyCurran65}
\begin{equation}
\label{eq:pde:general}
\frac{\partial}{\partial t} c_i - \adi D_i \Delta c_i
  = \sum_{j=1}^m N_{ij} \avj \r_j, \qquad
   x\in \Omega, \quad 0 < t < \tend.
\end{equation}
The parameter $D_i$ is the diffusion coefficient of species $Y_i$. It is 
considered homogeneous. For Cartesian coordinates the Laplacian
$\Delta$ is defined as
\[
\Delta = \frac{\partial^2}{\partial x_1^2} 
   +\frac{\partial^2}{\partial x_2^2} 
   +\frac{\partial^2}{\partial x_3^2}.
\]
The operator $\Delta$ has the dimension of $1/\textrm{length}^2$ $[\mathrm m^{-2}]$.
$D_i$ has the
dimension of surface per time $[\mathrm{m}^2/\mathrm{s}]$
and is assumed to be the same
throughout volume $\Omega$. In order to assess the importance of diffusion of
a substance for the behavior of the system, we shall modulate the
diffusion processes.  To this aim, $D_i$ is multiplied by the
dimensionless modulation parameter $\adi$, which has the
value of 1 in the reference state. 

In addition to Eq.~\pref{eq:pde:general},
which holds in the bulk, represented by $\Omega$, we need to
specify the behaviour at the boundary $\dOmega$ of the 
set  $\Omega$.
For a species $Y_i$ for which the membrane is a passive
barrier we assume a closed boundary condition,
\begin{equation}
\label{eq:bc_closed}
-\adi D_i \frac{\partial c_i}{\partial \nu} = 0 \qquad \text{on }\dOmega.
\end{equation}
The vector $\nu$ is the normal vector of unit length pointing out
of $\Omega$, and $\partial c_i/\partial \nu$ is the derivative
of $c_i$ in the direction of $\nu$. This makes the quantity
$-\adi D_i \partial c_i/\partial \nu$ the outward diffusive flux of 
species $Y_i$. 

If a species is transported across the boundary of $\Omega$,
then the boundary condition takes a similar form, 
but with non-zero right-hand side,
\begin{equation}
\label{eq:bc_open}
-\adi  D_i \frac{\partial c_i}{\partial \nu} 
   = \afi  f_i \qquad \text{on }\dOmega,
\end{equation}
where $f_i$ is the transport rate per unit surface area. $f_i$ can
be due to active transport or a (bio)chemical reaction at the
membrane surface
and may be a function of the
concentrations $c_1\ld c_n$.
There is at least 
one non-zero flux $\afi  f_i$ (one of which we conveniently 
assume to correspond
to $i=1$). Like other $\alpha$'s,  $\afi $ is a parameter
that we shall use to modulate the transport activity. 
This modulation is equivalent
to changes in the activity or concentration of membrane enzymes or carriers.
We shall refer to 
the `flux $J$ through the system' as the magnitude of the net export of $Y_1$,
integrated over the surface of $\Omega$ and normalized by the surface area:
\begin{equation}
\label{def:J}
J = \frac{\alpha_{f_1}}{\mod{\dOmega}} \int_{\dOmega} f_1 \, dS
\end{equation}
We normalize the total flow with respect to the corresponding surface area
to obtain a surface-averaged flux. The difference in role
of a reaction rate in
the bulk phase ($\r_j$) and a transport rate ($f_i$) is underlined by the
difference in the dimensionality---the former is given in volume concentration
per time unit (e.g., M/s), whereas the latter is expressed in surface
concentration (density) per time unit (e.g., $\mathrm{M\cdot m/s}$). 
We comprise Eq.~\pref{eq:bc_closed} into Eq.~\pref{eq:bc_open} 
by allowing $f_i$ to be zero.
Some reactions that involve cytoplasmic substances only occur
in the membrane or at the membrane surface. 
This is equivalent to the substrate of the reaction being
exported and the product of the reaction being imported, and will
be treated as such.

\subsection{Parameter Modulation and Control Analysis}

\begin{figure}[ht]
\centerline{\psfig{figure=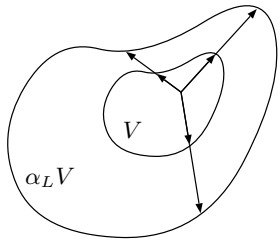,height=3.5cm}}
\caption{The domain $\Omega$ is scaled by geometric multiplication}
\label{fig:scaling_Omega}
\end{figure}

We have explicitly introduced
dimensionless modulation parameters $\avj$, $\adi $, and $\afi $, through which
the system can be modified. For the formulation and interpretation 
of control theorems associated with diffusion we need to
introduce two additional forms of modulation, associated with
space and time, in the following way. We assume that the spatial
and temporal variables $x_1,x_2,x_3$, and $t$, are related to
a reference spatial and temporal frame $\xi_1,\xi_2,\xi_3$, and $\tau$, 
via dimensionless modulation parameters $\al$ and $\at$:
\begin{equation}
\label{spacetime_modulation}
x_1= \al \xi_1, \qquad x_2 = \al \xi_2, \qquad x_3 = \al \xi_3, \qquad
\text{and} \qquad  t = \at \tau.
\end{equation}
Note that a reference time interval $[0,\Tend]$ is rescaled to
an actual time interval $[0,\at \Tend] = \at [0,\Tend]$; 
and similarly, the reference
volume $V$ is scaled to an actual volume $\al V$ 
(Figure~\ref{fig:scaling_Omega}).

We shall be explicit about the cause-and-effect relationship that
exists between \emph{parameters} and \emph{state}.
The parameters are the kinetic constants in rate 
functions $\r_j$ and $f_i$ and diffusion coefficients $D_i$,
together with their modulators $\avj$, $\adi $, and $\afi $. We also consider
the volume $\Omega$ and the time interval $[0,\Tend]$ as parameters, together
with the space and time modulators $\al$ and $\at$. 
A final `parameter' that we need to take into account is the
distribution of all species in space at the initial time.

For a given choice of all parameters the concentration of each of the species
$Y_i$ is determined uniquely as a function of space and time, according
to Eqs.~\pref{eq:pde:general} and~\pref{eq:bc_open}.
The choice of parameters yields a a unique set of
functions $c_i(x,t)$, each defined at every
point $x\in \al \Omega$ and at each time $t \in [0,\at\Tend]$. 

Note that the functions $c_i$ are defined on different domains
for different $\al$ and $\at$. For comparison purposes it will be
more convenient to think of $c_i$ as functions of the
reference coordinates $\xi$ and $\tau$ and the modulation parameters $\alpha$,
i.e.\ $c_i(\xi,\tau,\alpha)$,
where $\xi\in V$, $\tau \in [0,\Tend]$, and $\alpha = (\avj,\adi,\afi,\al,\at)$. 
The functions $c_i(\xi,\tau,\alpha)$ determine the \emph{state} of the system
associated with the choice of parameters. 

The change of coordinates, from $(x,t)$ to $(\xi,\tau)$, induces
a small modification in Eqs.~\pref{eq:pde:general}
and~\pref{eq:bc_open},
which now become
\begin{alignat}2
\label{eq:pde:general_XT}
&\at^{-1}\frac{\partial}{\partial \tau} c_i - \al^{-2} \adi  D_i \Delta c_i
  = \sum_{j=1}^m N_{ij} \avj \r_j, &\qquad
   &\xi\in \Omega, \quad 0 < \tau < \Tend \\
\label{eq:bc_open_XT}
&-\al^{-1}\adi  D_i \frac{\partial c_i}{\partial \nu} 
   = \afi  f_i &\qquad &\xi\in \dOmega, \quad 0 < \tau < \Tend
\end{alignat}
Note that the differentiation in the operators $\Delta$ and 
$\partial/\partial\nu$ is now with respect to~$\xi$.
The definition of the flux, however, remains unchanged under this
change of coordinates, because of the surface averaging:
\begin{equation}
\label{def:J_XT}
J = \frac{\alpha_{f_1}}{\mod{\dOmega}} \int_{\dOmega} f_1 \, dS
\end{equation}

In the vein of Metabolic Control Analysis 
we will consider a \emph{reference state} associated with the choice 
$\avj=\adi =\afi =\al=\at=1$ 
 and a second state, resulting from an infinitesimal 
change in one of these modulators. The difference between the 
two states will be taken to characterize 
the control that the modulated parameter has on the state of the
system. The dimensionless parameter $\avj$ is the tool we shall use
to modulate the activity of reaction $\r_j$.
Departures of $\avj$ from 1 correspond to a modulation of enzyme catalytic
activity or concentration, or to a proportional modulation of all rate
constants of the reaction step if it is not enzyme catalyzed (this
modulation leaves the equilibrium constant unchanged, for details 
see~\cite{KholodenkoWesterhoff93,KholodenkoWesterhoff94,KholodenkoWesterhoff95}, 
where this approach was introduced). For instance, if 
reaction~\pref{reaction:general}
is a protein-protein interaction (a central reaction in signaling
networks), the rate $\r_j$ depends on two kinetic constants, $k_{\mathrm{on}}$
and $k_{\mathrm{off}}$,
and a change in $\avj$ corresponds to equal proportional change of both
constants that leaves the dissociation (equilibrium) constant $K_d$
unchanged. The ratio of the fractional change in a state variable 
of the reaction network,
such as $J$ or $c_k$, and the fractional change in
$\avj$ determines the control coefficient with respect to reaction rate $\r_j$
(in the limit of infinitesimally small changes): 
\[
C_{v_j}^J = \left. \frac{d \ln J}{d\ln\avj}\right|_{\avj=1}, \qquad
C_{v_j}^{c_k} = \left.\frac{d \ln c_k}{d\ln\avj}\right|_{\avj=1},
\]
where $\avj$ = 1
corresponds to the reference state (all other parameters are assumed
fixed). Similarly, considering two states corresponding to small
modulations in dimensionless parameters, $\afi$ and $\adi$, we define the
control coefficients with respect to transport reaction and diffusion,
$C^J_{f_i}$ and $C^J_{D_i}$.  

An intrinsically novel control coefficient emerges as the
heterogeneous spatial organization of cellular reaction network is
taken into account. This coefficient is obtained when the size ($L$) of
the system is modified through a modulation of the dimensionless
parameter $\al$, as shown schematically in Fig.~\ref{fig:scaling_Omega}, 
whilst keeping fixed
all other parameters. In a cell setting, these parameters include the
volume and surface concentrations of enzymes and other molecular
forms, the total concentrations of which are conserved (for instance, the total
amount of phosphorylated and unphosphorylated forms of a
protein) in network
reactions described by Eqs.~\pref{eq:pde:general_XT} and~\pref{eq:bc_open_XT}.
Therefore, increasing the cell size also implies adding
conserved chemical moieties and additional enzymes both to the bulk
aqueous phase and cell membranes. Accordingly, the `spatial' control
coefficient is defined as $C^J_L = d\ln J/d\ln \al$, 
where only the spatial variables are
modulated, as described above (Eq.~\pref{spacetime_modulation}). 
It is of note, that the
parameter $L$ can be interpreted as the characteristic cell size,
e.g., the volume to surface ratio, whose fractional changes are equal
to the fractional changes in $\al$.

\medskip

A special case of Eq.~\pref{eq:pde:general_XT} 
arises when we only consider steady-state,
non-equilibrium configurations; these correspond to solutions 
$c_1(\xi)\ld c_n(\xi)$ of the steady-state equations
\begin{equation}
\label{eq:pde:stationary}
- \al^{-2} \adi  D_i \Delta c_i
  = \sum_{j=1}^m N_{ij} \avj \r_j, \qquad
   \xi\in \Omega,
\end{equation}
with boundary conditions~\pref{eq:bc_open_XT}.
For the purposes of this paper we assume that 
solutions of this set of equations are locally unique and asymptotically
stable.

\section{Results: the first summation theorem}

In a metabolic network where concentrations are homogeneous,
the sum of the flux control by all biochemical reactions
equals~$1$~\cite{KacserBurns73,HeinrichRapoport73}.
We generalize this summation theorem for systems where
concentrations may not be homogeneous and where transport occurs ---
here the sum over all the reactive processes need not be 1. 
The flux may be partially controlled by diffusion or by transport
and it is because of this control that the classical summation theorem
no longer holds.
The generalized summation theorem
takes the form
\begin{equation}
\label{th:1}
\sum_{i=1}^n C_{D_i}^J + \sum_{j=1}^m C_{v_j}^J +\sum_{i=1}^n C_{f_i}^J
  = 1.
\end{equation}
In words this theorem reads: the total control by all diffusion, reaction,
and transport processes on any steady state flux  
equals~$1$. The proof of Eq.~\pref{th:1} is based on Euler's
theorem on homogeneous function, which can be stated as follows.
Let $g$ be a function of $p_1\ld p_n$, such that for all $p_i$ and
for all $\lambda >0$, 
\begin{equation}
\label{eq:scaling_f}
g(\lambda^{\beta_1} p_1\ld \lambda^{\beta_n} p_n) 
= \lambda^\gamma g(p_1\ld p_n),
\end{equation}
for some $\beta_1,\ldots,\beta_n,\gamma\in\R$,
i.e.,\ $g(p_1,\ldots,p_n)$ is the same function
after the transformation
\begin{equation}
\label{eq:scaling_g}
\tilde g= \lambda^\gamma g,
  \qquad \tilde p_i = \lambda^{\beta_i} p_i.
\end{equation}

Then
\begin{equation}
\label{eq:lemma:scaling}
\sum_{i=1}^n \beta_i \, \frac{\partial \ln |g|}{\partial \ln p_i} = \gamma.
\end{equation}

It may be noted that for some of the parameters $p_i$, $\beta_i$ 
may equal zero.
These parameters are not modulated in Eq.~\pref{eq:scaling_g}, 
and Eq.~\pref{eq:lemma:scaling} shows that these parameters
are absent from the summation theorem.

For $\gamma$ and all $\beta_i$ equal to~$1$, 
Eq.~\pref{eq:lemma:scaling} reminds of the flux control
summation theorem of Metabolic Control 
Analysis~\cite{KacserBurns73,HeinrichRapoport73}, but it may
not be clear \emph{a priori} what function and which parameters should be
considered.  When searching for other than the traditional summation
theorems, a strategy may prove useful.  One strategy is that of
leaving the system in essence in the same state (cf.~\cite{KacserBurns73}).  
In terms of the system under consideration this translates into
the concentrations $c_i(\xi,\tau)$ remaining the same when
the parameters $p_i$ of Eq.~\pref{eq:scaling_g} are changed.

In order to implement the above strategy and examine under what type
of parameter changes $c_i(\xi,\tau)$ might be constant, one may inspect the
equations that define~$c_i$.  For steady states 
Eqs.~\pref{eq:pde:stationary} and~\pref{eq:bc_open_XT} define
$c_i(\xi)$.  Multiplying both equations by the same factor $\lambda$ should
leave their solutions unchanged.  Consequently, multiplying all $\adi$, all
$\avj $, and all $\afi $ by that same factor should not change the solutions
$c_i(\xi)$ either: the terms in the equations are homogeneous functions of first
order of the parameters $\adi$, $\avj $, and $\afi $.  
When only these parameters are
modulated, the specific transport rate remains unaltered, because the
metabolite concentrations remain constant and neither 
$\adi $, $\avj$, or $\afi $ occurs in the function $f_i(c_1\ld c_n)$.
Consequently $f_i$ remains unchanged, and therefore we can calculate the
corresponding flux as a function of $\lambda$,
\[
J(\lambda) = \lambda \, \frac{\alpha_{f_1}}{\mod{\dOmega}} 
    \int_{\dOmega} f_1 = \lambda J(1).
\]
We can reformulate this result as follows: we have shown that the flux
$J$, as a function of the parameters $\adi$, $\avj$, and $\afi$, satisfies
\begin{multline*}
J(\lambda \alpha_{D_1}\ld \lambda \alpha_{D_n}, 
  \lambda \alpha_{v_1}\ld \lambda \alpha_{v_m}, 
  \lambda \alpha_{f_1}\ld \lambda \alpha_{f_n}) \\
= \lambda J(\alpha_{D_1}\ld \alpha_{D_n}, \alpha_{v_1}\ld \alpha_{v_m}, \alpha_{f_1}\ld
  \alpha_{f_n}).
\end{multline*}
Equality~\pref{th:1} now follows from Eq.~\pref{eq:lemma:scaling}.

It should be remembered that the parameters $\adi $, $\avj$, and
$\afi $ represent the activities of diffusion, chemical, and transport
processes, respectively.  Consequently, Eq.~\pref{th:1} really refers to a
theorem concerning the control by diffusion, chemistry and transport.
It is an extension, both in concept and in proof, of the familiar flux
control summation theorem of Metabolic Control Analysis.  The latter
only dealt with the middle terms of the left-hand side and then only
the enzyme catalyzed reactions thereof, it assumed concentrations to
be homogeneous in space, and neglected explicit transport terms.  

A further extension to the above methodology will generalize the above
summation theorem to systems that depend on time: Dropping
the assumption of steady state, Eqs.~\pref{eq:pde:general_XT}
and~\pref{eq:bc_open_XT} become the equations
that define $c_i(\xi,\tau)$.  Multiplying both these equations by the same
factor $\lambda$ again leaves their solutions unchanged. Consequently:
\begin{equation}
\label{th:1a}
-C^{J(\tau)}_t + \sum_{i=1}^n C_{D_i}^{J(\tau)} 
  + \sum_{j=1}^m C_{v_j}^{J(\tau)} +\sum_{i=1}^n C_{f_i}^{J(\tau)}
  = 1.
\end{equation}
Here $C^{J(\tau)}_t$ coincides with the time-control coefficient defined 
in~\cite{AcerenzaSauroKacser89}. How can we understand this theorem and
in particular this new time-control coefficient?  Let us
first consider the situation that would ensue from simultaneously and
equally increasing all the parameters $\adi$, $\avj$, and
$\afi$ by $1\,\%$ in the
absence of an increase in $\at$.  Because all activities increase, also the
flux $J$ should increase.  However, the system should also change
more quickly in time, i.e.\ the increase in magnitude of $J$ should occur
earlier.  To obtain a proportionate increase in $J$, one should look at
an earlier point in time, i.e.\ time should be $1\, \%$ earlier,
or, keeping in mind that $t=\at T$, $\at = 0.99$.

$C^{J(\tau)}_t$ should not be confused with the time-dependent
control coefficients defined by \cite{WesterhoffVanDam87,HeinrichReder91}.
For a further discussion of the distinctions, see 
also~\cite{KholodenkoDeminWesterhoff97}
for the control analysis of relaxations
in the vicinity of the steady state.  

Often, as time goes to infinity,
the system relaxes to a steady state.  The control by time on the flux
then reduces to zero and the summation theorem for the time dependent
control coefficients reduces to the one for steady state control
coefficients (Eq.~\pref{th:1}).

\medskip
The summation theorem~\pref{th:1a} is more general than suggested by the
above derivation.  It can be derived without requiring explicit
equations for the time evolution of the system, such as 
Eqs.~\pref{eq:pde:general_XT} and~\pref{eq:bc_open_XT}.
The proof then considers an actual physical transformation of the
system under consideration: all the time dependent elemental processes
are increased by the same factor $\lambda$.  This implies that everything
happens in the same way, but faster by that same factor $\lambda$.  
We represent this concept of `faster' by the parameter $\at$,
as above (since $t = \at \tau$,  $\at<1$ corresponds
to `faster', and $\at>1$ to `slower').
Writing the parameter-dependence explicitly,
\[
c_i(\xi,\tau) =c_i(\xi,\tau;\adi,\avj,\afi,\al,\at),
\]
we have
\[
c_i\left(\xi,\tau;\lambda\adi,\lambda\avj ,\lambda\afi,\al,\frac\at\lambda\right)
 = c_i(\xi,\tau;\adi,\avj,\afi,\afi,\at),
\]
or in other words, $c_i$ is homogeneous of zeroth order in 
$\adi $, $\avj $, $\afi $, and $1/\at$.  
Using Eq.~\pref{eq:lemma:scaling}
this yields the summation theorem for concentrations,
\begin{equation}
\label{th:1b}
- C^{c_k}_t + \sum_{i=1}^n C_{D_i}^{c_k} 
  + \sum_{j=1}^m C_{v_j}^{c_k} +\sum_{i=1}^n C_{f_i}^{c_k}
  = 0.
\end{equation}
where $c_k$ is short for $c_k(\xi,\tau)$.
A similar result for the flux,
\[
J\left(\tau;\lambda\adi ,\lambda\avj ,\lambda\afi ,\frac\at\lambda\right)
 = \lambda J(\tau;\adi ,\avj ,\afi ,\at),
\]
shows that the flux is homogeneous of first order in 
$\adi $, $\avj $, $\afi $, and $1/\at$, resulting in the summation
theorem Eq.~\pref{th:1a}.

\section{A novel summation theorem}

The above theorem is an extension of the theorem 
of~\cite{AcerenzaSauroKacser89}
to cases with transport and diffusion.  The theorem we shall develop
now is intrinsically new as its origin lies in spatial aspects that
have hitherto been neglected in Metabolic Control Analysis~\cite{Fell97}.  
The
strategy we follow is again to consider a transformation for which the
system properties should not change.  We
first discuss the nature of diffusion, as illustrated by Einstein's
diffusion equation:
\[
\delta x = \sqrt{n_D D \delta t},
\]
where $n_D$ is a number of order 1 that depends on the dimensionality.
In words, the displacement due to diffusion increases with the square
root of the diffusion coefficient.  Consequently, increasing all
diffusion coefficients by $\lambda^2$ should increase all displacement by the
factor $\lambda$.  

This becomes obvious when considering Eq.~\pref{eq:pde:general_XT}.
Replacing $\adi$ by  $\lambda^2\adi$ and $\al$ by $\lambda\al$ leaves 
Eq.~\pref{eq:pde:general_XT} unchanged. The boundary 
condition~\pref{eq:bc_open_XT} in a natural way suggests
replacing $\afi $ by $\lambda\afi $, 
so that this equation is also essentially
unchanged. As a result this choice of parameters leaves
the solution $c_i$ invariant:
\begin{equation}
\label{eq:invariant_c_i}
c_i(\xi,\tau; \lambda^2\adi, \avj, \lambda\afi,\lambda \al, \at)
 = c_i(\xi,\tau; \adi, \avj, \afi, \al,\at).
\end{equation}
We immediately deduce the second summation theorem for concentrations,
\begin{equation}
\label{th:2_c}
2\sum_{i=1}^n C^c_{D_i} + \sum_{i=1}^n C^c_{f_i} + C^c_L = 0,
\end{equation}
where $c$ stands for any $c_i(\xi,\tau)$. Similarly, since
the flux satisfies
\[
J(\tau;\lambda^2\adi, \avj, \lambda\afi,\lambda \al, \at)
 = \lambda J(\tau;\adi, \avj, \afi, \al,\at),
\]
the associated summation theorem results
\begin{equation}
\label{th:2_j}
2\sum_{i=1}^n C^{J(\tau)}_{D_i} + \sum_{i=1}^n C^{J(\tau)}_{f_i} + C^{J(\tau)}_L = 1.
\end{equation}
Naturally, if we only consider steady states, then $c_i(\xi,\tau)$ reduces
to $c_i(\xi)$, and $J(\tau)$ to $J$.

The modulation of $V$, via $\al$, implicated
in the evaluation of the corresponding control coefficient is one in
which the size of the system changes; it is important,
however, that the system size changes whilst
keeping the volume concentration of enzymes in the bulk aqueous phase
constant and the surface concentrations of the membrane enzymes
constant.

\bigskip

\begin{figure}[ht]
\centerline{\psfig{figure=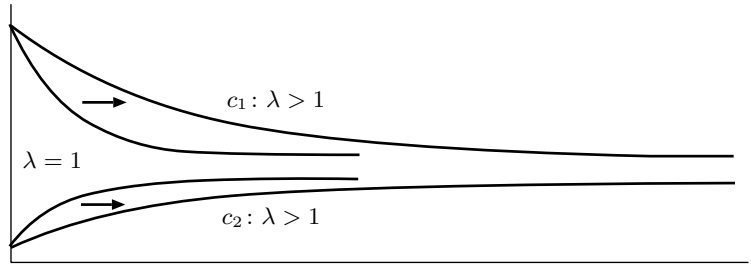,height=3.5cm}}
\caption{The scaling in $\lambda$ for a 
two-species system}
\label{fig:scaling_c}
\end{figure}

To illustrate the modulation introduced above, consider
the case of a one-dimensional spatial domain of the 
form $(0,\infty)$; the membrane is thought to be at $x=0$,
and the cell is deemed so large that the `other end' is effectively
at infinity. For such a `half-infinite' cell, in non-equilibrium 
steady-state, Eq.~\pref{th:2_j} reduces to
\begin{equation}
\label{th:2a}
2 \sum_{i=1}^n C^J_{D_i} + \sum_{i=1}^n C^J_{f_i} = 1
\end{equation}
where $J$ is the flux at the border ($x=0$); $C^J_L$ vanishes
as a result of the infinite size. Figure~\ref{fig:scaling_c}
shows how the solution $c_i(x)$ depends on the parameter modulation
by the factor $\lambda$ given in Eq.~\pref{eq:invariant_c_i}. Note that
here $c_i$ is plotted as a function of $x$; when plotted as a function
of $\xi$, the graphs for different~$\lambda$ coincide
(cf.~\pref{eq:invariant_c_i}).

\drop{
The choice of the modulators implies that functions $c_i(\xi)$ 
that are solutions of Eqs.~\pref{eq:pde:stationary} and~\pref{eq:bc_open_XT} 
for $\lambda=1$ remain solutions when $\lambda>1$. 
By substituting into Eq.~\pref{def:J_XT} we find that
\[
J(\lambda) = \lambda\, J(1),
\]
so that Eq.~\pref{th:2a} follows from Eq.~\pref{eq:lemma:scaling}.

}

\section{An explicit example system}
\label{sec:ex}

For certain simple cases control coefficients can be calculated explicitly.
We do so here for a kinase/phosphatase example 
that we also discussed in~\cite{BrownKholodenko99,KholodenkoBrownHoek00}.
There are two species, $\X$ and $\XP$, which are the unphosphorylated
and phosphorylated form of a given protein; the kinase converts
$\X$ into $\XP$, and the phosphatase does the reverse.

\begin{figure}
\centerline{\psfig{figure=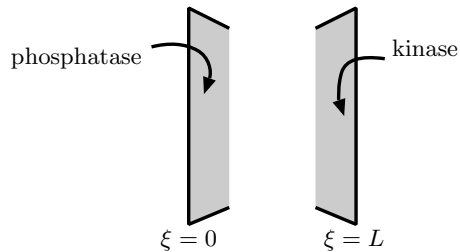,height=3cm}}
\caption{Geometry 1: a slice bounded by parallel membranes}
\label{fig:kp1}
\end{figure}

\begin{figure}
\centerline{\psfig{figure=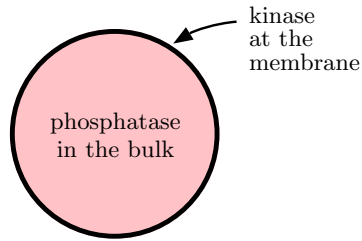,height=3cm}}
\caption{Geometry 2: a spherical cell}
\label{fig:kp2}
\end{figure}

Besides providing an illustration of the concepts and results 
discussed above we wish to demonstrate the influence
of the system geometry on the control coefficients. We consider
two cases:
\begin{enumerate}
\item $\Omega$ is bounded by two parallel membranes; the kinase 
reaction is localized to one membrane, the phosphatase to the other
(Figure~\ref{fig:kp1});
\item $\Omega$ is a spherical region (e.g.\ a cell or an organelle)
bounded by a membrane; the kinase is localized at the membrane,
but the phosphatase is distributed throughout the bulk 
(Figure~\ref{fig:kp2}).
\end{enumerate}

As we shall see below, in the case of parallel membranes the
size control coefficient $C^J_L$ is \emph{negative}: an increase in the 
distance between the membranes reduces the flux through the system.
However, when the phosphatase is not membrane-bound but distributed,
as in the second case, increasing the system size
actually increases the flux: $C^J_L$ is \emph{positive}.

\subsection{Two parallel membranes}

We assume that the system is large with respect to the
distance $L$ between the membranes and therefore 
adopt a one-dimensional formulation. We also restrict the study 
to (non-equilibrium) steady state.

Since no reaction takes place in the bulk, the equations
satisfied by $\X$ and $\XP$ in the bulk (Eq.~\pref{eq:pde:stationary})
reduce to
\begin{equation}
\label{eq:DE_int}
-\al^{-2}\ad D\, \XP'' = -\al^{-2}\ad D\, \X'' = 0, \qquad 0<\xi<L.
\end{equation}
Note that the spatial variable is $\xi$, as discussed above;
the physical variable $x$ is given by $x=\al\xi$.

For the kinase and phosphatase reactions we assume the rate
functions 
\[
\ak k_k (\kk\X-\XP) \qquad \text{and}\qquad \ap k_p (\XP-\kp\X).
\]
Note the dimensionality of $k_{k,p}$: the rate functions are
surface fluxes, but $\XP$ and $\X$ are bulk concentrations;
$k_k$ and $k_p$ therefore have dimension length/time.
The coefficients $\alpha_{k,p}$ and $\kappa_{k,p}$ are dimensionless.

The kinase and phosphatase reactions both enter the
description as boundary conditions (cf.~\pref{eq:bc_open_XT}):
\begin{alignat}2
 \begin{aligned}
  \al^{-1}\ad D\,\XP' &= \ap k_p (\XP-\kp\X) \\
  \al^{-1}\ad D\,\X' &= - \ap k_p (\XP-\kp\X)
 \end{aligned}
 & \qquad && \text{at $\xi=0$ (phosphatase)}
\label{eq:bc_0}
\\[5pt]
 \begin{aligned}
  -\al^{-1}\ad D\,\XP' &= -\ak k_k (\kk\X-\XP) \\
  -\al^{-1}\ad D\,\X' &=   \ak k_k (\kk\X-\XP)
 \end{aligned}
 & \qquad && \text{at $\xi=L$ (kinase)}
\label{eq:bc_L}
\end{alignat}

Equation~\pref{eq:DE_int} implies that the spatial gradients
of $\XP$ and $\X$ (i.e.~$\XP'$ and~$\X'$) are constant in $\xi$;
let us set $\delta = \XP'$.
From the solution of equations~(\ref{eq:DE_int}--\ref{eq:bc_L})
the relevant information for our purposes is $\delta$:
\begin{equation}
\label{def:delta}
\delta = \frac{M(\beta_p\ak\kk-\beta_k\ap\kp)}
              {\al^{-1}\ad D(\beta_k+\beta_p) + L\beta_k\beta_p}.
\end{equation}
Here $M$ is the total concentration $\X+\XP$, which is independent of
$\xi$ by Eq.~\pref{eq:DE_int}, and 
\[
\beta_k = \ak  k_k(1+\kk) \qquad \text{and} \qquad
\beta_p = \ap  k_p(1+\kp).
\]

Because of the spatial separation of the reactions
the flux as defined by Eq.~\pref{def:J} equals the diffusive
flux,
\begin{equation}
\label{eq:flux_delta}
J = -\al^{-1} \ad D\, \XP' = -\al^{-1}\ad D\, \delta.
\end{equation}
From Eqs.~\pref{def:delta} and~\pref{eq:flux_delta} we can determine
the control coefficients of the flux with respect to 
the kinase and phosphatase reactions and to diffusion, as well as the
`length'-control coefficient:%
\footnote{Note that 
these formulas contain no $\alpha$'s, since they result from evaluation at 
$\al=\ad=\ak=\ap=1$; consequently, $\beta_{k,p} = k_{k,p}(1+\kappa_{k,p})$.}
\begin{align*}
C^J_k &= \frac{D\beta_p}{D(\beta_k+\beta_p) + L\beta_k\beta_p},\\
C^J_p &= \frac{D\beta_k}{D(\beta_k+\beta_p) + L\beta_k\beta_p},\\
C^J_D &= \frac{L\beta_k\beta_p}{D(\beta_k+\beta_p) + L\beta_k\beta_p},\\
C^J_L &= -\frac{L\beta_k\beta_p}{D(\beta_k+\beta_p) + L\beta_k\beta_p}.
\end{align*}
Clearly these expressions satisfy the
summation theorems Eqs.~\pref{th:1} and~\pref{th:2_j}:
\begin{align*}
&C^J_D + C^J_k + C^J_p = 1,\\
&2C^J_D + C^J_k + C^J_p + C^J_L = 1.
\end{align*}
Note that both the kinase and the phosphatase
reaction are boundary effects, and therefore give rise to 
boundary control coefficients $C^J_{f_i}$ of Eqs.~\pref{th:1} and~\pref{th:2_j}.

For concentration control we choose to consider not a pointwise
concentration, but the difference between concentrations at
opposite ends:
\[
\XP(L)- \XP(0) = L\delta.
\]
Denoting the control coefficients with respect to this
quantity as $C^c$, we have
\begin{align*}
C^c_k &= \frac{D\beta_p}{D(\beta_k+\beta_p) + L\beta_k\beta_p},\\
C^c_p &= \frac{D\beta_k}{D(\beta_k+\beta_p) + L\beta_k\beta_p},\\
C^c_D &= -\frac{D(\beta_k+\beta_p)}{D(\beta_k+\beta_p) + L\beta_k\beta_p},\\
C^c_L &= \frac{D(\beta_k+\beta_p)}{D(\beta_k+\beta_p) + L\beta_k\beta_p}.
\end{align*}

\begin{figure}[ht]
\centering
\centerline{\psfig{figure=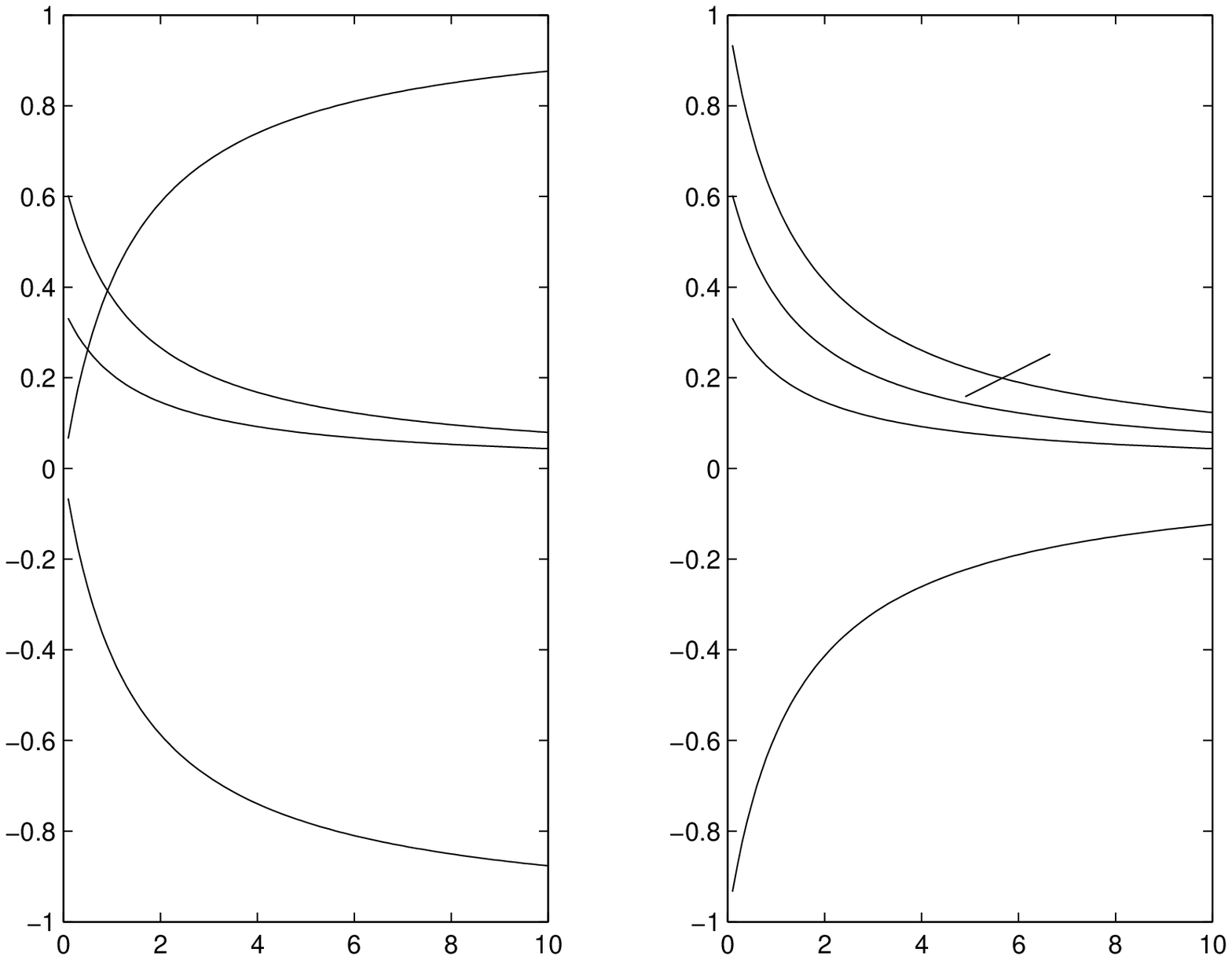,height=7cm}}
\vskip5mm
\caption{Values of control coefficients as a function
of $L$ for the first example (parallel membranes at distance
$L$, the kinase and phosphatase are localized at different membranes;
$D = 1\,\mu m^2/s$, $k_p = k_k = 1/s$, $\kp = 0.1$, $\kk = 10$).}
\label{fig:ex_cc1}
\end{figure}

An example of the values of these control coefficients is given in
Figure~\ref{fig:ex_cc1}.

\subsection{A spherical cell}
For the second example we consider the same system as
above, but now in a spherical cell of radius $L$.
We assume spherical symmetry throughout.

The kinase reaction is again localized at the membrane,
and therefore results in a boundary condition
\begin{equation}
\begin{aligned}
  -\al^{-1}\ad D\,\XP' &= \ak k_k (\XP-\kk\X), \\
  -\al^{-1}\ad D\,\X' &= - \ak k_k (\XP-\kk\X),
 \end{aligned}
 \qquad  \text{at $\xi=L$ (kinase).}
\label{eq:bc_L2}
\end{equation}
The spatial variable is again $\xi$, but it now represents a \emph{radial}
coordinate.
The phosphatase reaction, which in this example is distributed 
throughout the bulk, now enters as a right-hand side in the differential
equation for $\X$ and $\XP$ (cf.~\pref{eq:pde:stationary}):
\[
\begin{aligned}
&{-\al^{-2}\ad D}\, \xi^{-2}(\xi^2\XP')' =  -\ap k_p (\XP - \kp \X),\\
&{-\al^{-2}\ad D}\, \xi^{-2}(\xi^2\X')' =  \ap k_p (\XP - \kp \X),
\end{aligned}
\qquad 0< \xi< L.
\]
Again note the dimensions of $k_k$ and $k_p$: similar to the
previous example, $k_k$ has dimension length/time; but here $k_p$ converts
a bulk concentration into a bulk flux, and therefore $k_p$ has
the (more familiar) dimension 1/time.
As before, the total concentration is independent of $\xi$,
and denoted by $M$:
\begin{equation}
\label{eq:M}
M= \XP+\X
\end{equation}

As a first step we investigate different boundary conditions, where we 
prescribe the flux $J$:
\begin{equation}
\label{bc:prescribed_flux} 
-\al^{-1}\ad D\, \XP' = J = \al^{-1}\ad D\, \X',
\end{equation}
at $\xi=L$. For this choice of boundary conditions we can solve the
system exactly. Set $u=\XP-\kp \X$; then $u$ satisfies the equation
\[
-\al^{-2}\ad D\, \xi^{-2}(\xi^2 u')' = -\ap k_p (1+\kp) u,
\]
with boundary condition
\[
-\al^{-1}\ad D\,  u'(L) = (1+\kp)J.
\]
We find that $u$ is given by
\begin{equation}
\label{eq:sol_u}
u(\xi) = -J\, \frac {1+\kp}{\al^{-1}\ad D}     \; 
   \frac {L^2}{\gamma L \cosh \gamma L - \sinh \gamma L}
    \, \frac {\sinh \gamma \xi}\xi,
\end{equation}
where the inverse length scale $\gamma$ is defined by
\[
\gamma^2 = \al^2 \, \frac{\ap k_p(1+\kp)}{\ad D}.
\]
It follows that
\begin{equation}
\label{def:c_i}
\begin{aligned}
\XP(\xi) &=  \frac{\kp}{1+ \kp }\, M
         + \frac 1{1+\kp} \, u(\xi)\\
\X(\xi) &=  \frac1{1+\kp }\, M
         - \frac 1{1+\kp} \, u(\xi).
\end{aligned}
\end{equation}

To return to the boundary 
conditions~\pref{eq:bc_L2} we equate the
flux $J$ in Eq.~\pref{bc:prescribed_flux} with the right-hand side 
of Eq.~\pref{eq:bc_L2}:
\[
J = \ak k_k (\XP(L) - \kk \X(L)) = 
   \ak k_k \left (\frac{1+\kk}{1+\kp}\,u(L) 
                + \frac{\kp-\kk}{1+\kp} \, M\right);
\]
using the expression~\pref{eq:sol_u} for $u$ this yields
\begin{equation}
\label{eq:sol_J}
J = \frac{\nu}{1+\beta} ,\qquad
  \nu = \ak k_k M \, \frac{\kp-\kk}{1+\kp}, \qquad 
  \beta = \ak k_k \,\frac{1+ \kk}{\al^{-1}\ad D} \;
\frac L {\gamma L \coth \gamma L - 1}
\end{equation}
($\coth x$ is the hyperbolic cotangent,
$\cosh x/\sinh x$).
\drop{
It follows from Eq.~\pref{eq:sol_J}
that
\[
1 + \beta = \frac {\nu M} {\XP(L) - \kk \X(L)}.
\]
}

As in the previous example we thus obtain explicit
formulas for the control of the various parameters on the flux:
\begin{align*}
C_{k}^J &= \frac1{1+\beta} \\
C_{D}^J &= \frac12 \, \frac{\beta}{1+\beta} \, 
  \frac{ \gamma L \coth \gamma L -2 + \gamma^2 L^2(\coth^2\gamma L-1)}%
                              {\gamma L\coth \gamma L - 1} \\
C_p^J &= \frac12 \, \frac{\beta}{1+\beta} 
          \, \frac{\gamma L \coth \gamma L - \gamma^2 L^2(\coth^2\gamma L-1)}%
                  {\gamma L \coth \gamma L - 1} \\
C_L^J &= \frac{\beta}{1+\beta} 
          \, \frac{1 - \gamma^2 L^2(\coth^2\gamma L-1)}%
                  {\gamma L \coth \gamma L - 1} \\
\end{align*}

As an example of the concentration control values we choose
a point at distance $L/2$ from the center of the cell. Denoting
the control coefficients with respect to this quantity as
$C^c$, we calculate 
\begin{align*}
C_{k}^c &= \mu C^J_k \\
C_{D}^c &= \mu \left(C^J_D + \frac12 \phi -1\right)\\
C_p^c &= \mu \left(C^J_p - \frac12 \phi \right)\\
C_L^c &= \mu \left(C^J_L - \phi +1 \right)
\end{align*}
Here
\[
\mu = \frac{u(L/2)}{\kp M + u(L/2)} \qquad \text{and} 
\qquad \phi = \frac{\gamma^2 L^2}{\gamma L \coth \gamma L - 1}
    -\frac{\gamma L}2 \coth \left(\frac{\gamma L}2\right).
\]

\begin{figure}[ht]
\centering
\centerline{\psfig{figure=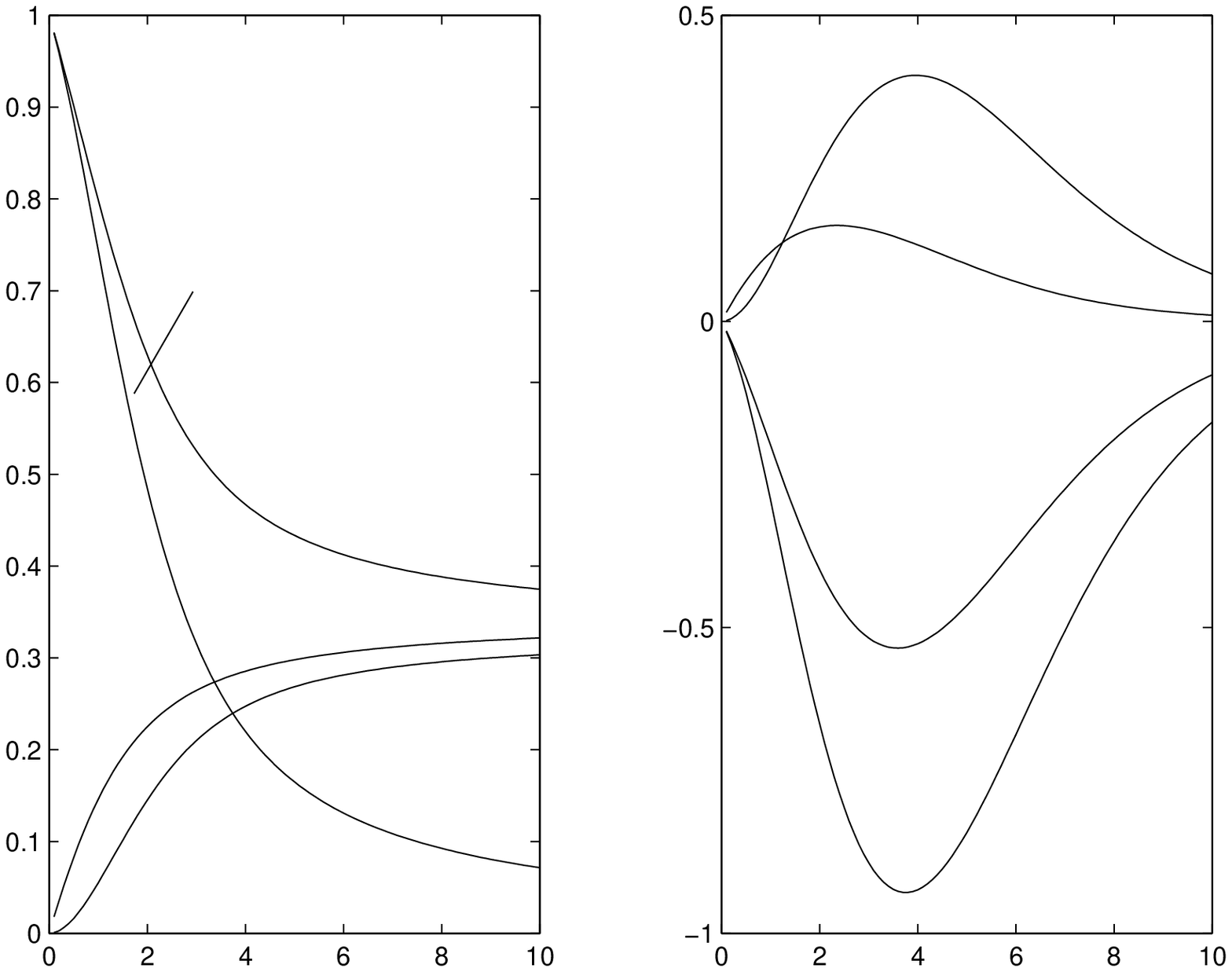,height=7cm}}
\vskip5mm
\caption{Values of control coefficients as a function
of $L$ for the second example (a spherical cell with a kinase on the membrane
and a phosphatase in the cytosol;
$D = 1\,\mu m^2/s$, $k_p = k_k = 1/s$, $\kp = 0.1$, $\kk = 10$).
$C^c$ denotes control of the concentration measured halfway between the center and the membrane. }
\label{fig:ex_cc2}
\end{figure}

An example of the values of these control coefficients is given
in Figure~\ref{fig:ex_cc2}.

\section{Discussion}

Until fairly recently, biochemistry focused more on the time dependent
aspects of processes than on the spatial aspects.  Exceptions were the
membrane-mediated subcompartmentation of cellular metabolism.  The
cellular subcompartments were mostly considered to be homogeneous in
terms of concentrations and the heterogeneity between them was
analyzed in terms of the activity of transport catalysts in the
membranes~\cite{WesterhoffVanDam87}.  
Inhomogeneity of metabolites and even ions within
aqueous subcompartments has often been proposed~\cite{Kell79}
but calculations
have shown that at least for central metabolic routes aqueous
diffusion should be fast on the time scale of the catalytic turnover
of the enzymes also given the proximity of the enzyme molecules due to
the small sizes of cells~\cite{WesterhoffVanDam87,WesterhoffWelch92}.
Metabolic Control Analysis was
developed from this perspective and proved a useful way to rationalize
the 
study of regulation and control of intermediary metabolism~\cite{Fell97} and
free-energy transduction~\cite{WesterhoffVanDam87}.  
More than an analysis procedure,
Metabolic Control Analysis also provided biochemists with a set of
laws (theorems) that govern control of fluxes and concentrations in
metabolic networks.  One type of these were the summation theorems,
stating that the sum of the control by all the individual enzyme
catalyzed reactions on flux and
concentration~\cite{KacserBurns73,HeinrichRapoport73}, oscillation
amplitude~\cite{KholodenkoDeminWesterhoff97} or relaxation time 
constant~\cite{HeinrichReder91} 
should equal 1, 0, 1, and 1 respectively.

Those were the days where biochemistry focused on metabolism.  However, only a
minor fraction of the known genomes encodes primary metabolism.  Much
coding capacity is devoted to regulation, in part through signal
transduction~\cite{Kell79,KholodenkoWesterhoffSchwaberCascante00}.
Signal transduction has (i) the spatial aspect
that many signals arrive at the cell's plasma membrane and have to be
transferred to its nuclear DNA, (ii) the added complication that much
of the signal is transferred by proteins in direct phosphoryl
transfers, rather than with small molecules as messengers, and (iii)
the feature that those proteins are often present at a thousand times
lower concentrations than the usual metabolite.  All these three
features work in the direction that might make diffusion of the
components limiting in signal transfer.  Indeed, it was calculated
that for realistic parameter
values, concentration gradients in signal 
(defined as in~\cite{FranckePostmaWesterhoffBlomPeletierTA}) should
arise~\cite{Fell80,BrownKholodenko99,KholodenkoBrownHoek00}, 
which might even force the cell to take refuge to alternative
mechanisms of movement of the signal proteins~\cite{Kholodenko02}.

With the realization of the possible importance of transport and
diffusion for cell biochemistry, there came a need to update Metabolic
Control Analysis so as to include the aspects of diffusion.  This is
what the first part of this paper accomplished for the summation
theorem.  The result was rather simple, to the sum of the enzyme control
coefficients also coefficients for the control by the diffusion of the
various species needed to be added, as well as control coefficients for
the transport processes and for the non-enzyme-catalyzed chemical
reactions, before it amounted to 1 for flux control and 0 for
concentration control.  This result is perhaps not surprising, but
has been overlooked effectively in the literature for quite a while.
In studies of the control of growth rate of \emph{E. coli} by the lactose
permease, 
it was found that a sum of 1 could not be found and it
was suspected that part of the control resided in diffusion across the
outer membrane~\cite{JensenWesterhoffMichelsen93,DykhuizenDeanHartl87}.
It was found that none of the glycolytic
enzymes had significant control on the glycolytic flux in yeast.  This
was considered a problem until strong indications arose of
substantial control in glucose transport into the 
cells~\cite{ReijengaSnoepDiderichVanVerseveldWesterhoffTeusink01}.
Thus, if anything, our extension of the summation theorem to include
the transport steps may open eyes to the possibility that
transport controls metabolic fluxes.  Similarly we should expect
renewed attention to arise from the other term in our summation
theorem, i.e.\ the one that refers to control by the diffusion.  Of
course the novelty lies not so much in the possibility of such 
control, but in the
fact that it may account for a shortfall of metabolic control from
the total of~1.

We have also formulated the summation theorem for the case of
inhomogeneity in both time and space, i.e.\ where a spatial gradient
may be developing over time.  In that case diffusion and time
turn out to share in the total flux control of~1.  In the case of
homogeneity in time (i.e.\ at steady state), the former contribution
disappears, in case of spatial homogeneity the latter.  The
experimental system that springs to mind is that of
Belousov-Zhabotinsky reactions exhibiting time varying 
spirals~\cite{ZhabotinskyZaikin73}, but
closer to biology the wave type oscillations in yeast extracts~\cite{MairWarnkeMuller01}.
The most relevant application may be in developmental biology, but for
this an extension to Hierarchical Control 
Analysis~\cite{SnoepVanderWeijdenAndersenWesterhoffJensen02} will be
needed.

\bigskip
While the first flux summation result, Eq.~\pref{th:1}, and
the first concentration summation result, Eq.~\pref{th:1b}, are
a simple extension of existing theory, the summation theorems
concerning changes of spatial scale (Eqs.~\pref{th:2_c},
\pref{th:2_j}, and~\pref{th:2a}) are fundamentally new,
and introduce a concept that is new to Metabolic Control Analysis:
control by size.

\medskip

`Control by the size', the quantity measured by $C_L$, 
is a concept with some interesting properties. 
For instance, $C^J_L$ and $C^c_L$ can be both positive and negative.
This is related to the fact that the size of a system can
influence the `productive' capacity of the system in
different ways at the same time. Let us examine this for flux control:
\begin{itemize}
\item Size can be an obstacle: if material has to travel over a distance $L$, 
as in the example of parallel
membranes (Case 1 of Section~\ref{sec:ex}), 
then an increase in~$L$ implies a decrease of the (average) concentration gradient,
and therefore causes a reduction
in diffusive flux. This is the dominant effect in the large-$L$ limit
of Case 1, as shown in Fig.~\ref{fig:ex_cc1}, where it leads to a limit value
of $C^J_L$ of $-1$.
\item On the other hand, size can also be
a resource: if a reaction takes place in a bulk region
of linear size $L$, then $L$ is also a measure of
reactive capacity: an increase in size will result in
an increase of total reactive flux. This is demonstrated by the second
example, where the input- and output boundaries are co-located and
no material is forced to travel over a distance $L$.
Here the control by size is positive, and in the limit of small
$L$ the control even approaches $1$.
\end{itemize}

The new summation theorems have a number of interesting implications.
Perhaps the most striking one is that they impose an upper limit to
the control that diffusion may have on the flux through the system.  The
first summation theorem suggested that this number be 1, i.e.\ if the
biochemical and the transport processes were in excess and the system
were large; the second, however, shows that it is only $1/2$ 
whenever the control by size is positive (i.e.\ when increasing the size
results in a higer flux).  Clearly, in such cases there must always 
be another process
controlling flux, in addition to diffusion.  
In the example of Fig.~\ref{fig:ex_cc2}
the other processes were the biochemical reactions.  

The example of Figure~\ref{fig:scaling_c} also shows how
changing the diffusion rate alters the concentration profiles. 
By increasing the diffusion rate we increase the width of a
`boundary layer', in which
a significant part of the reactions takes place. As a result
the flux is also increased. Note, however, that the relative increase in
size of the boundary layer (which scales as $\lambda$, since 
$x = \al \xi = \lambda \xi$) is the square root of the increase
of the diffusion coefficients ($\adi = \lambda^2$).
This can be explained intuitively as follows. 
The size of the boundary layer is
determined by counteracting forces; a thin boundary layer
provides for fast spatial transport (on the scale of the layer
thickness), but little reactive
capacity, while for a thick layer the situation is inversed.
An increase in diffusion rate is therefore used partially
to increase transport rate, and partially to increase
reaction capacity. The number $1/2$ is the result of
this trade-off.

\noindent\textsc{Acknowledgement}
\par\penalty1000
\noindent 
This work was supported by the RTN network 
`Nonlinear Partial Differential Equations describing Front 
Propagation and other Singular Phenomena',
HPRN-CT-2002-00274, and by the National Institute of Health through Grant
GM59570.

\newpage
\bibliographystyle{bphj}
\bibliography{ref}
\end{document}